\documentclass[12pt]{article}
\usepackage[margin = 20mm]{geometry}
\usepackage{amsmath}
\usepackage{amsfonts}
\usepackage{mathrsfs}
\usepackage{amssymb}
\usepackage{amsthm}
\usepackage{url}
\usepackage{caption}
\usepackage{cases}
\usepackage{blkarray}
\usepackage{enumerate}
\usepackage{tikz}
\usepackage{tkz-graph}
\usetikzlibrary{shapes.geometric}
\usetikzlibrary{babel}
\usepackage{pgfplots}
\pgfplotsset{compat=1.5}
\usepackage{tocloft}
\usepackage{tabularx}
\usepackage{setspace}
\usepackage{array,color,colortbl}   

\usepackage{abstract}
\usepackage{nicematrix}
\usepackage{xcolor, colortbl}
\usepackage{float}
\usepackage[ruled]{algorithm2e}
\LinesNumbered
\newtheorem{thm}{Theorem}[section]

\newtheorem{lem}[thm]{Lemma}

\theoremstyle{definition}

\newtheoremstyle{break}
{\topsep}{\topsep}%
{}{}%
{\bfseries}{}%
{\newline}{}%
\theoremstyle{break}
\numberwithin{equation}{section}
\def\eref#1{$(\ref{#1})$}
\def\sref#1{\S$\ref{#1}$}
\def\lref#1{Lemma~$\ref{#1}$}
\def\tref#1{Theorem~$\ref{#1}$}

\def\fref#1{Figure~$\ref{#1}$}

\def\tbref#1{Table~$\ref{#1}$}
\def\spacer{{\vrule height 2.25ex width 0ex depth1.0ex}}
\def\nauty{\FuncSty{Nauty}}

\newcommand{\jn}{\Vert}

\renewcommand{\geq}{\geqslant}
\renewcommand{\leq}{\leqslant}
\renewcommand{\ge}{\geqslant}
\renewcommand{\le}{\leqslant}

\renewcommand{\epsilon}{\varepsilon}
\newcommand{\Z}{\mathbb{Z}}
\newcommand{\E}{\mathcal{E}}
\newcommand{\F}{\mathcal{F}}
\renewcommand{\P}{\mathcal{P}}
\renewcommand{\L}{\mathcal{L}}
\newcommand{\LL}{\mathscr{L}}

\newcommand{\T}{\mathcal{T}}
\newcommand{\M}{\mathcal{M}}
\newcommand{\po}{\prec}
\newcommand{\poe}{\preccurlyeq}

\makeatletter
\g@addto@macro\bfseries{\boldmath}
\makeatother

\title{Perfect $1$-factorisations of $K_{11,11}$}
\author{Jack Allsop\\
	\small Institut für Mathematik\\[-0.5ex]
	\small Freie Universität Berlin\\[-0.5ex]
	\small Berlin 14195, Germany\\
	\small\tt allsop@mi.fu-berlin.de
	\and
	Ian M. Wanless\\
	\small School of Mathematics\\[-0.5ex]
	\small Monash University\\[-0.5ex]
	\small Vic 3800, Australia\\
	\small\tt ian.wanless@monash.edu}
\date{}

\begin{document}
	
	\maketitle
	
	\begin{abstract}
		A perfect $1$-factorisation of a graph is a decomposition of that graph into $1$-factors such that the union of any two $1$-factors is a Hamiltonian cycle. A Latin square of order $n$ is row-Hamiltonian if for every pair $(r,s)$ of distinct rows, the permutation mapping $r$ to $s$ has a single cycle of length $n$.  We report the results of a computer enumeration of the perfect $1$-factorisations of the complete bipartite graph $K_{11,11}$.  This also allows us to find all row-Hamiltonian Latin squares of order $11$. Finally, we plug a gap in the literature regarding how many row-Hamiltonian Latin squares are associated with the classical families of perfect $1$-factorisations of complete graphs.
	\end{abstract}
	
	\section{Introduction}\label{s:intro}
	
	A \emph{$1$-factor}, or \emph{perfect matching}, of a graph $G$ is a
	set of edges of $G$ with the property that every vertex of $G$ is in
	exactly one of the edges. A \emph{$1$-factorisation} of $G$ is a
	partition of its edge set into $1$-factors. Let $\F$ be a
	$1$-factorisation of $G$ and let $f$ and $f'$ be distinct $1$-factors
	in $\F$. The edges in $f$ and $f'$ together form a subgraph of
	$G$ which is a union of cycles of even length. If $f \cup f'$ induces
	a Hamiltonian cycle in $G$, regardless of the choice of $f$ and $f'$,
	then $\F$ is a \emph{perfect $1$-factorisation}. Two
	$1$-factorisations $\F$ and $\E$ of $G$ are \emph{isomorphic} if there
	exists a permutation $\phi$ of the vertices of $G$ which maps the set
	of $1$-factors in $\F$ onto the set of $1$-factors in $\E$. In this
	case, $\phi$ is an \emph{isomorphism} from $\F$ to $\E$. An
	\emph{automorphism} of $\F$ is an isomorphism from $\F$ to itself. The
	\emph{automorphism group} of $\F$ is the set of all automorphisms of
	$\F$ under composition.
	
	The main purpose of this paper is to report the results of a computer
	enumeration of the perfect $1$-factorisations of the complete
	bipartite graph $K_{11,11}$. It is known that a perfect
	$1$-factorisation of $K_{n,n}$ can only exist if $n=2$ or $n$ is odd
	(see, e.g.,~\cite{rowhamcount}). It is conjectured that a perfect
	$1$-factorisation of $K_{n,n}$ does exist in these cases. However,
	this conjecture is a long way from being resolved. There are few known
	infinite families of perfect $1$-factorisations of complete bipartite
	graphs~\cite{nu4,newatom,p1fbip}, and these only cover graphs
	$K_{n,n}$ where $n\in\{p,2p-1,p^2\}$ for some odd prime $p$.  Up to
	isomorphism there are $1$, $1$, $1$, $2$ and $37$ perfect
	$1$-factorisations of $K_{2,2}$, $K_{3,3}$, $K_{5,5}$, $K_{7,7}$ and
	$K_{9,9}$, respectively~\cite{rowhamcount}.
	
	Perfect $1$-factorisations of complete bipartite graphs are
	related to perfect $1$-factorisations of complete graphs
	(see~\cite{symm1facs} for full details of this relationship).
	In particular, a perfect $1$-factorisation of $K_{2n}$ can be
	used to build a perfect $1$-factorisation of $K_{2n-1,2n-1}$
	using a construction that we call Kotzig's construction, which
	is given explicitly in \sref{s:GA}.  So the existence of a
	perfect $1$-factorisation of $K_{2n}$ implies the existence of
	a perfect $1$-factorisation of $K_{2n-1, 2n-1}$, but it is not
	known whether the converse holds.  In 1964,
	Kotzig~\cite{Kotzconj} famously conjectured that a perfect
	$1$-factorisation of $K_{2n}$ exists for all positive integers
	$n$. This conjecture remains even further from resolution than
	the conjecture on the existence of perfect $1$-factorisations
	of complete bipartite graphs. There are only three known
	infinite families of perfect $1$-factorisations of complete
	graphs~\cite{newatom}, and these only cover graphs $K_{2n}$
	where $2n\in \{p+1, 2p\}$ for an odd prime $p$.  Up to
	isomorphism there are $1$, $1$, $1$, $1$, $1$, $5$, $23$ and
	$3155$ perfect $1$-factorisations of $K_2$, $K_4$, $K_6$,
	$K_8$, $K_{10}$, $K_{12}$, $K_{14}$ and $K_{16}$,
	respectively~\cite{P1F14,P1F10,P1F16,P1F162,P1F12}.
	
	The main result of this paper is the following theorem.
	
	\begin{thm}\label{t:main}
		There are $687\,121$ perfect $1$-factorisations of
		$K_{11,11}$ up to isomorphism. Of these, $2657$ have a
		non-trivial automorphism group.
	\end{thm}

	The structure of this paper is as follows. In \sref{s:alg} we
	discuss our enumeration algorithm for proving
	\tref{t:main}. There is an equivalence between
	$1$-factorisations of complete bipartite graphs and Latin
	squares. As a result, the catalogue behind \tref{t:main}
	allows us to enumerate several interesting classes of Latin
	squares of order $11$, as discussed in \sref{s:LS}. In
	\sref{s:invar}, we discuss how useful various invariants are
	for distinguishing our enumerated objects. In \sref{s:GA}, we
	prove a new property of a well known family of perfect
	$1$-factorisations of complete graphs.
	
	To reduce the risk of programming errors, all computations
	described in this paper were performed independently by each
	author, then crosschecked.  The combined computation time was
	under two CPU years.
	
	\section{The algorithm}\label{s:alg}
	
	In this section we describe how we generated the perfect
	$1$-factorisations of $K_{11,11}$. The algorithm we used is similar
	to the algorithm used in \cite{P1F16} to generate the perfect
	$1$-factorisations of $K_{16}$.
	
	A \emph{partial $1$-factorisation} of a graph $G$ is a
	collection of pairwise disjoint $1$-factors of $G$. Let $\P$
	be a partial $1$-factorisation of $G$ and let $f$ and $f'$ be
	distinct $1$-factors in $\P$. If $f \cup f'$ induces a
	Hamiltonian cycle in $G$ then $(f, f')$ is a \emph{perfect
		pair}. If every pair of distinct $1$-factors in $\P$ is
	perfect, then $\P$ is called perfect. An \emph{ordered}
	partial $1$-factorisation is a partial $1$-factorisation with
	an order on its $1$-factors.  We use $\F=[f_1,f_2,\dots,f_a]$
	to denote an ordered partial $1$-factorisation with
	$1$-factors $f_1,\dots,f_a$ and then $\F\jn f_{a+1}$ to denote
	$[f_1,f_2,\dots,f_a,f_{a+1}]$, the ordered partial
	$1$-factorisation obtained by appending $f_{a+1}$ to $\F$. Two
	ordered partial $1$-factorisations $\F=[f_1,f_2,\ldots,f_a]$
	and $\E=[e_1,e_2,\ldots,e_a]$ of $G$ are \emph{isomorphic} if there
	is a permutation $\psi$ of $\{1,2,\ldots,a\}$ and a
	permutation $\phi$ of the vertices of $G$ which maps $f_i$
	onto $e_{\psi(i)}$ for every $i\in\{1,2,\ldots,a\}$.
	In this section, we will primarily be discussing ordered
	$1$-factorisations of $K_{n,n}$. However, for most of our
	purposes the order will be inconsequential, so we will often
	just refer to such objects as $1$-factorisations. 
	
	Throughout this section, we will assume that the vertices of
	$K_{n,n}$ are labelled by $u_1,u_2,\ldots,u_n$ and
	$v_1,v_2\ldots,v_n$, where there is an edge between $u_i$ and
	$v_j$ for all $\{i, j\}\subseteq\{1,\dots,n\}$.  For brevity
	we will write the edge $\{u_i,v_j\}$ as $u_iv_j$, and
	similarly for other graphs. We will call an isomorphism
	\emph{direct} if it preserves $\{u_1,u_2,\ldots,u_n\}$ and
	$\{v_1,v_2\ldots,v_n\}$ setwise, and \emph{indirect} if it
	exchanges these two sets.
	
	We now define a partial order $\poe$ on the set of ordered
	partial $1$-factorisations of $K_{n,n}$.  Let
	$\F=[f_1,f_2,\ldots,f_a]$ and $\E=[e_1,e_2,\ldots,e_b]$ be two
	distinct such partial $1$-factorisations. If $f_i=e_i$ for all
	$i\le\min(a,b)$ then $\F$ and $\E$ are
	incomparable. Otherwise, let $j$ be minimal such that $f_j
	\neq e_j$.  If $j>4$ then we deem $\F$ and $\E$ incomparable.
	So suppose that $j\le4$.  The edges in $f_j$ can be written as
	$u_1x_1,u_2x_2,\ldots,u_nx_{n}$ where
	$\{x_1,\dots,x_n\}=\{v_1,\ldots,v_n\}$.  Similarly the edges
	in $e_j$ can be written as $u_1y_1,u_2y_2,\ldots,u_ny_n$ where
	$\{y_1,\dots,y_n\}=\{v_1,\ldots,v_n\}$.  Let $\ell$ be minimal
	such that $x_\ell\neq y_\ell$.  We say that $\F \po \E$ if
	$x_\ell<y_\ell$ (in the lexicographical ordering). If
	$x_\ell>y_\ell$, we say that $\E\po\F$.  Let $\poe$ denote the
	reflexive closure of $\po$.
	
	Let $\F=[f_1,f_2,\ldots,f_a]$ be an ordered partial
	$1$-factorisation of $K_{n,n}$ with $a\geq4$. Denote by $\F^i$
	the ordered partial $1$-factorisation $[f_1,f_2,\dots,f_i]$. Say that
	$\F$ is \emph{minimal} if $\F^4\poe\E^4$ for every ordered partial
	$1$-factorisation $\E$ of $K_{n,n}$ that is isomorphic to
	$\F$. Note that if $\F=[f_1,f_2,\ldots,f_a]$ is a minimal perfect
	partial $1$-factorisation then
	\begin{equation}\label{e:f1f2}
		\begin{aligned}
			&f_1 = \{u_1v_1,u_2v_2,\dots,u_{n-1}v_{n-1},u_nv_n\} \text{ and,} \\
			&f_2 = \{u_1v_2,u_2v_3,\dots,u_{n-1}v_n,u_nv_1\}.
		\end{aligned}
	\end{equation}
	
	A graph isomorphism between vertex coloured graphs is
	\emph{colour preserving} if the colour of each vertex matches
	the colour of its image under the isomorphism.
	The software \nauty~\cite{nauty} is a practical algorithm for
	testing whether there is a colour preserving graph isomorphism
	between two vertex coloured graphs.  Isomorphism testing for
	$1$-factorisations of bipartite graphs can be converted into
	an isomorphism problem on vertex coloured graphs as
	follows. For a $1$-factorisation $\F=[f_1,f_2,\ldots,f_a]$ of
	a graph $G\subseteq K_{n,n}$ we construct a coloured graph
	$C(\F)$ containing
	\begin{itemize}
		\item green vertices $f_1,f_2,\ldots,f_a$ each joined to a blue vertex $F$,
		\item green vertices $u_1,\dots,u_n$ each joined to a red vertex $U$,
		\item green vertices $v_1,\dots,v_n$ each joined to a red vertex $V$,
		\item one black vertex for each edge in $G$ which is joined to one
		green vertex in each of the previous three categories to indicate
		the end points of the edge and the $1$-factor that contains the
		edge.
	\end{itemize}
	It is routine to check that two partial $1$-factorisations $\F$ and
	$\E$ are isomorphic if and only if there is a colour preserving graph
	isomorphism from $C(\F)$ to $C(\E)$. Also, the automorphism group of
	$\F$ is (group) isomorphic to the group of colour preserving
	automorphisms of $C(\F)$, which \nauty\ counts.   As an aside, the whole
	construction can be varied in an obvious way to solve the isomorphism
	problem for $1$-factorisations of non-bipartite graphs.
	
	Our algorithm for generating the perfect $1$-factorisations of
	$K_{n,n}$ is described in Procedure 2, and its subroutine
	\FuncSty{AddFactor} described in Procedure 1.  Steps 2 and 7
	of Procedure 2 can be handled in a straightforward manner
	using \nauty\ as discussed above, and represent a negligible
	fraction of the computation time.  As mentioned at the
	beginning of this section, our algorithm is similar to the one
	used in~\cite{P1F16} to generate the perfect
	$1$-factorisations of $K_{16}$. Apart from obvious adaptations
	to the bipartite setting, the main change is a refinement on
	when minimality checks are performed.

	\begin{algorithm}[h]
		\DontPrintSemicolon
		\SetAlgorithmName{Procedure}{}{}
		\SetKwInOut{Input}{input}
		\SetKwFunction{proc}{AddFactor$(n,\P,\T)$}
		\Input{An odd integer $n\ge5$\\
			A perfect partial $1$-factorisation $\P$ of $K_{n,n}$\\
			A set $\T$ of $1$-factors $t$ for which $\P\jn t$
			is a perfect partial $1$-factorisation
		}
		\SetKwProg{myproc}{Procedure}{}{end}
		\myproc{\proc}{
			\eIf{$|\P|=n$}
			{
				Output $\P$\;
			}
			{
				Let $e$ be an edge of $K_{n,n}\setminus\bigcup\P$ that is in the fewest $1$-factors in $\T$\;
				\For{$t \in \T$ containing $e$}
				{
					Let $\T^*$ be the set of $1$-factors $t^* \in \T$ such that $(t,t^*)$ is a perfect pair\;
					$\textnormal{\FuncSty{AddFactor}($n$, $\P\jn t$, $\T^*$)}$\;
				} 
			}
		}
		\caption{Recursively add $1$-factors to a perfect partial $1$-factorisation}
	\end{algorithm}
	
	\begin{algorithm}[h]
		\DontPrintSemicolon
		\SetAlgorithmName{Procedure}{}{}
		\SetKwInOut{Input}{input}
		\SetKwFunction{proc}{GenP1Fs$(n)$}
		\Input{An odd integer $n\ge5$}
		\SetKwProg{myproc}{Procedure}{}{end}
		\myproc{\proc}{
			Generate the set
			$\mathcal{S}$
			of minimal perfect partial $1$-factorisations of
			$K_{n,n}$ containing four $1$-factors and edges
			$u_1v_1$, $u_1v_2$, $u_1v_3$, $u_1v_4$\;
			\For{$P \in \mathcal{S}$}{
				Let $\T=\{1$-factors $t$ such that $P\jn t$ is a minimal perfect partial $1$-factorisation$\}$\;
				$\textnormal{\FuncSty{AddFactor}($n$, $P$, $\T$)}$\;
			}
			Screen the 1-factorisations output by \FuncSty{AddFactor} for isomorphism
		}
		\caption{Generate perfect $1$-factorisations of $K_{n,n}$}
	\end{algorithm}
	
	Our algorithm starts by producing the set $\mathcal{S}$ of all
	minimal perfect partial $1$-factorisations of $K_{n,n}$ that
	contain four $1$-factors and that include the edges $u_1v_1$,
	$u_1v_2$, $u_1v_3$ and $u_1v_4$. Note that $\poe$ is a total order
	on $\mathcal{S}$ and it follows that no two elements of $\mathcal{S}$
	are isomorphic. For each element
	$P\in\mathcal{S}$ we then find the set $\mathcal{T}$ of
	$1$-factors whose addition to $P$ preserves minimality and
	perfection. These $1$-factors are then used to recursively
	extend our partial $1$-factorisation by one $1$-factor at a
	time until it has been completed to a perfect
	$1$-factorisation. Other choices for the initial number of
	$1$-factors could have been possible. Our choice of four
	initial $1$-factors was made so that the cardinalities of both
	the sets $\mathcal{S}$ and $\mathcal{T}$ were manageable. For
	$n=11$ we found that $|\mathcal{S}|=13\,727\,482$ and
	$|\mathcal{T}|$ ranged from $13\,954$ down to 0. As a result of
	the minimality requirements, $|\mathcal{T}|$ generally trended
	downwards as $P$ increased in the $\poe$ order, resulting in
	significant speed-up for later parts of the computation.
	
	It is worth remarking that there are perfect partial
	$1$-factorisations of $K_{n,n}$ that consist of four
	$1$-factors (including the edges $u_1v_1$, $u_1v_2$, $u_1v_3$ and
	$u_1v_4$) but are not isomorphic to any element of
	$\mathcal{S}$. For example, form $\F_*$ from
	the two factors in \eref{e:f1f2} together with
	\begin{align*}
		f_3&=\{u_1v_3,u_2v_1,u_3v_5,u_4v_7,u_5v_4,u_6v_{10},u_7v_{11},u_8v_6,u_9v_2,u_{10}v_9,u_{11}v_8\},\text{ and}\\
		f_4&=\{u_1v_4,u_2v_{10},u_3v_1,u_4v_{11},u_5v_7,u_6v_9,u_7v_6,u_8v_2,u_9v_8,u_{10}v_3,u_{11}v_5\}.
	\end{align*}
	Note that $\F_*$ is perfect, and contains the edges
	$u_1v_1$, $u_1v_2$, $u_1v_3$ and $u_1v_4$. However, the minimal member
	of the isomorphism class of $\F_*$ is formed from 
	the two factors in \eref{e:f1f2} together with
	\begin{align*}
		f_3&=\{u_1v_3,u_2v_1,u_3v_5,u_4v_2,u_5v_9,u_6v_4,u_7v_6,u_8v_{11},u_9v_8,u_{10}v_7,u_{11}v_{10}\}
		,\text{ and}\\
		f_4&=\{u_1v_{10},u_2v_5,u_3v_7,u_4v_{11},u_5v_4,u_6v_3,u_7v_9,u_8v_6,u_9v_1,u_{10}v_2,u_{11}v_8\},
	\end{align*}
	and does not contain $u_1v_4$. So, the isomorphism class of
	$\F_*$ has no representative in $\mathcal{S}$.
	Notwithstanding this example, we next show that our algorithm
	performs the desired enumeration.
	
	\begin{lem}\label{l:algworks}
		The set of $1$-factorisations returned by
		\FuncSty{GenP1Fs}$(n)$ contains at least one
		representative from each isomorphism class of perfect
		$1$-factorisations of $K_{n,n}$.
	\end{lem}
	
	\begin{proof}
		Let $\M$ be an isomorphism class of ordered perfect
		$1$-factorisations of $K_{n,n}$. Since $\M$ is finite
		it must contain a minimal element
		$\F=[f_1,\dots,f_n]$. Minimality of $\F$ forces $f_i$
		to contain the edge $u_1v_i$ for $1\le i\le 4$
		(cf.~\eref{e:f1f2}), and hence
		$\F^4\in\mathcal{S}$. Let $U=\{f_5,f_6,\dots,f_n\}$.
		The minimality of $\F$ implies that $\F^4\jn f$ is
		minimal for $f\in U$. So $U\subseteq\mathcal{T}$ when
		\FuncSty{AddFactor} is called with input $\P=\F^4$.
		By induction on $k\in\{5,\ldots,n\}$, subsequent
		recursive calls to \FuncSty{AddFactor} will be made
		with input $\P$ that consists of $\F^4$ together with
		$k-4$ of the $1$-factors in $U$, whilst the
		$1$-factors in $U\setminus \P$ are in $\mathcal{T}$.
		In each inductive step the
		$1$-factor $t$ that is added to $\P$ will be whichever
		$1$-factor in $U$ contains the edge $e$ defined by
		Line 5 of \FuncSty{AddFactor}.  There must be such a
		$1$-factor available, because $U$ contains a
		$1$-factorisation of $K_{n,n}\setminus\F^4$, and
		$\mathcal{T}$ inherits a $1$-factorisation of the
		graph induced by whichever edges have not yet been
		included in $\P$.
		
		In the case $k=n$, we see that \FuncSty{AddFactor} will output an
		ordered factorisation that equals $\F$, up to the order of its
		$1$-factors. The result follows.
	\end{proof}	
	
	Both of our implementations of \FuncSty{GenP1Fs} were used to
	generate the perfect $1$-factorisations of $K_{n,n}$ for
	$n\in\{5,7,9,11\}$.  Results of both programs agreed with each
	other, and for $n\le9$ agreed with previously computed values
	\cite{rowhamcount}.  \tbref{tb:autP1F} shows the $687\,121$
	perfect $1$-factorisations of $K_{11,11}$ categorised by the
	size of their automorphism group. The third column of the
	table lists the number of direct automorphisms of each
	$1$-factorisation, and the fourth column lists the total
	number of automorphisms. The first column gives the number of
	isomorphism classes that attain the attributes listed in the
	row in question and the second column gives the number of such
	isomorphism classes that can be obtained from perfect
	$1$-factorisations of $K_{12}$ via Kotzig's construction.  The
	number of direct automorphisms of a $1$-factorisation $\F$ can
	be counted by using \nauty\ as described above, simply by
	changing the colour of the vertex $V$ to yellow so that it can
	no longer be interchanged with $U$ in any colour preserving
	automorphism of $C(\F)$.

	\begin{table}
		\hspace{0mm}
		\vspace{-5mm}
		\begin{center}
			\begin{tabular}[h]{|c|c|c|c|}
				\hline
				\spacer
				Count&From $K_{12}$&Direct automorphisms&Automorphisms\\
				\hline\spacer
				684464  &  0&   1 &     1\\
				100  &  15&   1 &     2\\
				2531  &  0&   2 &     2\\
				6  &   0&  5 &     5\\
				5  &   3&  5 &    10\\
				7  &  0&  10 &    10\\
				3  &  3&  10 &    20\\
				1  &  0&  22 &    22\\
				2  &  0&  55 &    55\\
				1  &  1&  55 &   110\\
				1  &  1&1210 &  2420\\
				\hline	
			\end{tabular}
			\caption{\label{tb:autP1F}Symmetries of perfect $1$-factorisations of $K_{11,11}$}
		\end{center}
	\end{table}

	\section{Latin squares}\label{s:LS}
	
	We start this section by giving some basic definitions regarding
	Latin squares, and explain their relationship to
	$1$-factorisations of complete bipartite graphs. We then apply
	this previously known theory to our new catalogue of perfect
	$1$-factorisations of $K_{11,11}$, uncovering some
	interesting Latin squares of order $11$.
	
	Let $n$ and $m$ be positive integers with $m \leq n$. An $m
	\times n$ \emph{Latin rectangle} is an $m \times n$ matrix of
	$n$ symbols, each of which occurs exactly once in each row and
	at most once in each column. A \emph{Latin square} of order
	$n$ is an $n \times n$ Latin rectangle. In this paper we will
	always assume that the rows and columns of a Latin square are
	indexed by its symbol set. Let $L$ be an $m \times n$ Latin
	rectangle. A \emph{subrectangle} of $L$ is a submatrix of $L$
	that is itself a Latin rectangle. A $k \times \ell$
	subrectangle is \emph{proper} if $1 < k \leq \ell < n$. A
	\emph{subsquare} of $L$ is a subrectangle of $L$ that is
	itself a Latin square.  A \emph{row cycle of length $k$} in
	$L$ is a $2 \times k$ subrectangle of $L$ that has no proper
	subrectangles. A \emph{row-Hamiltonian} Latin square is a
	Latin square that has no proper subrectangles. Equivalently, a
	Latin square of order $n$ is row-Hamiltonian if all of its row
	cycles have length $n$.  A related but strictly weaker
	property is named $N_\infty$, which applies to Latin squares
	that have no proper subsquares. Such properties are very
	natural for mathematicians to consider, but turn out to be
	quite elusive.  After more than 50 years of studying the
	existence question it has only very recently been established
	in \cite{Ninf} that $N_\infty$ Latin squares exist for all
	orders $n\notin\{4,6\}$.  The analogous but harder question
	for row-Hamiltonian Latin squares is very far from being
	solved, and provides one of the motivations for compiling
	catalogues for small orders.
	
	Let $L$ and $L'$ be Latin squares. If $L$ can be obtained from $L'$ by
	applying a permutation $\alpha$ to its rows, a permutation $\beta$ to
	its columns and a permutation $\gamma$ to its symbols, then $L$ and
	$L'$ are \emph{isotopic}, and $(\alpha, \beta, \gamma)$ is an
	\emph{isotopism} from $L'$ to $L$. Isotopism is an equivalence
	relation and the equivalence classes are called \emph{isotopism classes}.
	Latin squares in the same isotopism class have the same
	number of subrectangles of each dimension, so the row-Hamiltonian
	property is an isotopism class invariant.  An \emph{autotopism} of $L$
	is an isotopism from $L$ to itself. The \emph{autotopism group} of $L$
	is the set of all autotopisms of $L$ under composition.
	
	Let $L$ be a Latin square of order $n$. We can consider $L$ as
	a set of $n^2$ triples of the form
	$(\text{row},\text{column},\text{symbol})$, called
	\emph{entries}. A \emph{conjugate} of $L$ is a Latin square
	obtained from $L$ by choosing a permutation of $\{1, 2, 3\}$
	and applying it to the coordinates of each entry in
	$L$. Every conjugate of $L$ can thus be labelled by a $1$-line
	permutation of $\{1, 2, 3\}$, which gives the order of the
	coordinates of the conjugate relative to the order of the
	coordinates of $L$. For example, the $(213)$-conjugate of $L$
	is its matrix transpose. The $(132)$-conjugate of $L$ is its
	\emph{row-inverse}. If $L$ is isotopic to some conjugate of
	$L'$ then $L$ and $L'$ are \emph{paratopic}. A
	\emph{paratopism} from $L'$ to $L$ is a pair
	$(\mathcal{C},\mathcal{A})$ where $\mathcal{A}$ is a $1$-line
	permutation of $\{1,2,3\}$ specifying a conjugate $L''$ of
	$L'$, and $\mathcal{C}$ is an isotopism from $L''$ to
	$L$. Paratopism is an equivalence relation and the equivalence
	classes are called \emph{species}.  An \emph{autoparatopism}
	of $L$ is a paratopism from $L$ to itself. The
	\emph{autoparatopism group} of $L$ is the set of all
	autoparatopisms of $L$ under composition.
	
	Let $L$ be a Latin square. Let $\nu(L)$ be the number of
	conjugates of $L$ that are row-Hamiltonian. We will also say
	that $L$ has $\nu = \nu(L)$. Since the row-Hamiltonian
	property is an isotopism class invariant, it follows that
	$\nu$ is a species invariant. So if $\nu(L)=c$ then we will
	say that the species of Latin squares containing $L$ has
	$\nu=c$. Latin squares with $\nu=6$ are called
	\emph{atomic}. It is known~\cite{rowhamcount} that $\nu(L) \in
	\{0, 2, 4, 6\}$, since a Latin square is row-Hamiltonian if
	and only if its row-inverse is row-Hamiltonian.

	There is a natural equivalence between Latin squares of order
	$n$ and ordered $1$-factorisations of $K_{n,n}$. This
	equivalence is studied in \cite{rowhamcount,symm1facs}, for
	example, where the following observations are spelt out in
	detail.  Let $L$ be a Latin square of order $n$. Label the
	vertices in one part of $K_{n,n}$ by $\{c_1,c_2,\ldots,c_n\}$,
	corresponding to the columns of $L$, and the vertices in the
	other part by $\{s_1,s_2,\ldots,s_n\}$, corresponding to the
	symbols of $L$. For row $i$ of $L$, we define a $1$-factor
	$f_i$ of $K_{n,n}$ by adding the edge $c_js_k$ to $f_i$
	whenever $L_{i,j}=k$. Then $\F=[f_1,f_2,\ldots,f_n]$ is an
	ordered $1$-factorisation of $K_{n,n}$, where the order on the
	$1$-factors comes from the order of the rows of $L$. It is
	easy to see that this construction is also reversible, giving
	a map $\F\mapsto\L(\F)$ from ordered $1$-factorisations of
	$K_{n,n}$ to Latin squares of order $n$.  The subgraph of
	$K_{n,n}$ induced by the $1$-factors $f_i$ and $f_j$ is a
	union of cycles of even length, and it contains a cycle of
	length $2k$ if and only if there is a row cycle in $\L(\F)$ of
	length $k$ hitting rows $i$ and $j$. Thus $\F$ is perfect if
	and only if $\L(\F)$ is row-Hamiltonian.

	Let $\F$ and $\E$ be ordered $1$-factorisations of
	$K_{n,n}$. From the definition of $\L(\F)$ and $\L(\E)$, it is not hard to see that $\F$ is isomorphic to $\E$ if and only if $\L(\F)$ is
	isotopic to $\L(\E)$ or the row-inverse of $\L(\E)$ (see~\cite{symm1facs} for details). 	
	
	\begin{lem}\label{l:symmnu}
		Suppose that $L$ is any Latin square of order $n$ that is isotopic
		to its transpose. For $X\in\{123,132,213,231,312,321\}$ let $\F_X$
		denote the ordered $1$-factorisation of $K_{n,n}$ for which
		$\L(\F_X)$ is the $(X)$-conjugate of $L$. Then
		$\F_{123},\F_{132},\F_{213}$ and $\F_{231}$ are all isomorphic.
		Hence, if $L$ is row-Hamiltonian then $\nu(L)\in\{4,6\}$ and if the
		$(321)$-conjugate of $L$ is row-Hamiltonian then $\nu(L)\in\{2,6\}$.
	\end{lem}
	
	\begin{proof}
		The proof is similar to that of \cite[Lem.~5]{rowhamcount}.
		Any Latin square $L$ would have an indirect isomorphism from
		$\F_{123}$ to $\F_{132}$ and also from $\F_{213}$ to $\F_{231}$.
		The fact that $L$ is isotopic to its transpose means that $\F_{123}$
		is isomorphic to $\F_{213}$. Hence, $\F_{123},\F_{132},\F_{213}$ and
		$\F_{231}$ are isomorphic to each other. 
		Thus the following four statements are equivalent:
		\begin{itemize}
			\item $L$ is row-Hamiltonian,
			\item The row-inverse of $L$ is row-Hamiltonian,
			\item The transpose of $L$ is row-Hamiltonian,
			\item The $(231)$-conjugate of $L$ is row-Hamiltonian.
		\end{itemize}
		The lemma follows.
	\end{proof}

	\tbref{tb:smallnLSrowhamatomic} gives the number of species and
	isotopism classes containing row-Hamiltonian Latin squares, as well as
	the number of species containing atomic Latin squares of small
	orders. The data for orders up to 9 was determined by
	Wanless~\cite{rowhamcount}, and the number of species containing
	atomic Latin squares of order 11 was determined by Maenhaut and
	Wanless~\cite{atom11}. The number of species containing
	row-Hamiltonian Latin squares of order $n$ exactly matches the numbers
	of perfect $1$-factorisations up to isomorphism of $K_{n,n}$ for all
	$n \in \{2, 3, 5, 7, 9\}$. However, this trend does not continue for
	order $11$ as is shown concretely by \eref{e:fr4} below. It was
	observed in \cite{rowhamcount} that there are no Latin squares of
	order $n$ with $\nu=4$ for $n \leq 9$, and that this trend also does
	not continue for $n=11$.
	
	\begin{table}
		\hspace{0mm}
		\vspace{-5mm}
		\begin{center}
			\begin{tabular}[h]{|c|c|c|c|}
				\hline
				\spacer
				order&row-Hamiltonian species & row-Hamiltonian isotopism classes & atomic species \\
				\hline\spacer
				2&1&1&1\\
				3&1&1&1\\
				5&1&1&1\\
				7&2&2&1\\
				9&37&64&0\\
				11&687\,115&1\,374\,132&7\\
				\hline	
			\end{tabular}
			\caption{\label{tb:smallnLSrowhamatomic}Row-Hamiltonian and atomic Latin squares of small order}
		\end{center}
	\end{table}
	
	From the set of representatives of isomorphism classes of perfect
	$1$-factorisations of $K_{11,11}$, it is a simple task to obtain
	representatives of each species of row-Hamiltonian Latin squares of
	order $11$.  This can be achieved by using \nauty\ as described in
	\sref{s:alg}, except that we recolour the vertex $F$ red.
	Each autoparatopism group is also automatically calculated
	by \nauty, which allows us to deduce the following data. 
	
	\begin{thm}\label{t:LSenum}\mbox{ }
		\begin{itemize}
			\item There are $687\,115$ species containing row-Hamiltonian Latin squares of order $11$. Of these, $2660$ have a non-trivial autoparatopism group, $687\,096$ have $\nu=2$, $12$ have $\nu=4$ and $7$ have $\nu=6$.
			\item There are $1\,374\,132$ isotopism classes containing row-Hamiltonian Latin squares of order $11$. Of these, $5104$ have a non-trivial autotopism group.
		\end{itemize}
	\end{thm}

	\begin{table}
		\hspace{0mm}
		\vspace{-5mm}
		\begin{center}
			\begin{tabular}[h]{|c|c|c|c|c|}
				\hline
				\spacer
				Count&From $K_{12}$&Autotopisms&Autoparatopisms&$\nu$\\
				\hline\spacer
				684455&0&1&1&2\\
				99&14&1&2&2\\
				8&0&1&2&4\\
				1&1&1&2&6\\
				2531&0&2&2&2\\
				5&0&5&5&2\\
				4&3&5&10&2\\
				1&0&5&10&6\\
				1&0&10&10&2\\
				1&0&10&10&4\\
				2&0&10&20&4\\
				2&2&10&20&6\\
				1&0&22&22&2\\
				1&1&10&60&6\\
				1&0&55&110&4\\
				1&1&55&110&6\\
				1&1&1210&7260&6\\
				\hline	
			\end{tabular}
			\caption{\label{tb:autpargrps}Symmetries of species of row-Hamiltonian Latin squares of order $11$}
		\end{center}
	\end{table}

	\tref{t:LSenum} allows us to fill in two previously unknown
	entries in the last row of
	\tbref{tb:smallnLSrowhamatomic}. \tbref{tb:autpargrps} shows
	the $687\,115$ species of row-Hamiltonian Latin squares of
	order $11$ classified according to how much symmetry they
	have. In that table the second column lists how many species
	contain a symmetric Latin square whose $(321)$-conjugate is
	row-Hamiltonian, indicating that it can be obtained from one
	of the perfect $1$-factorisations of $K_{12}$ via Kotzig's
	construction.  The third and fourth columns give the orders of
	the autotopism group and the autoparatopism group,
	respectively. The last column gives the value of $\nu$ and the
	first column reports how many species attain the attributes
	listed in the row in question.
	
\begin{figure}
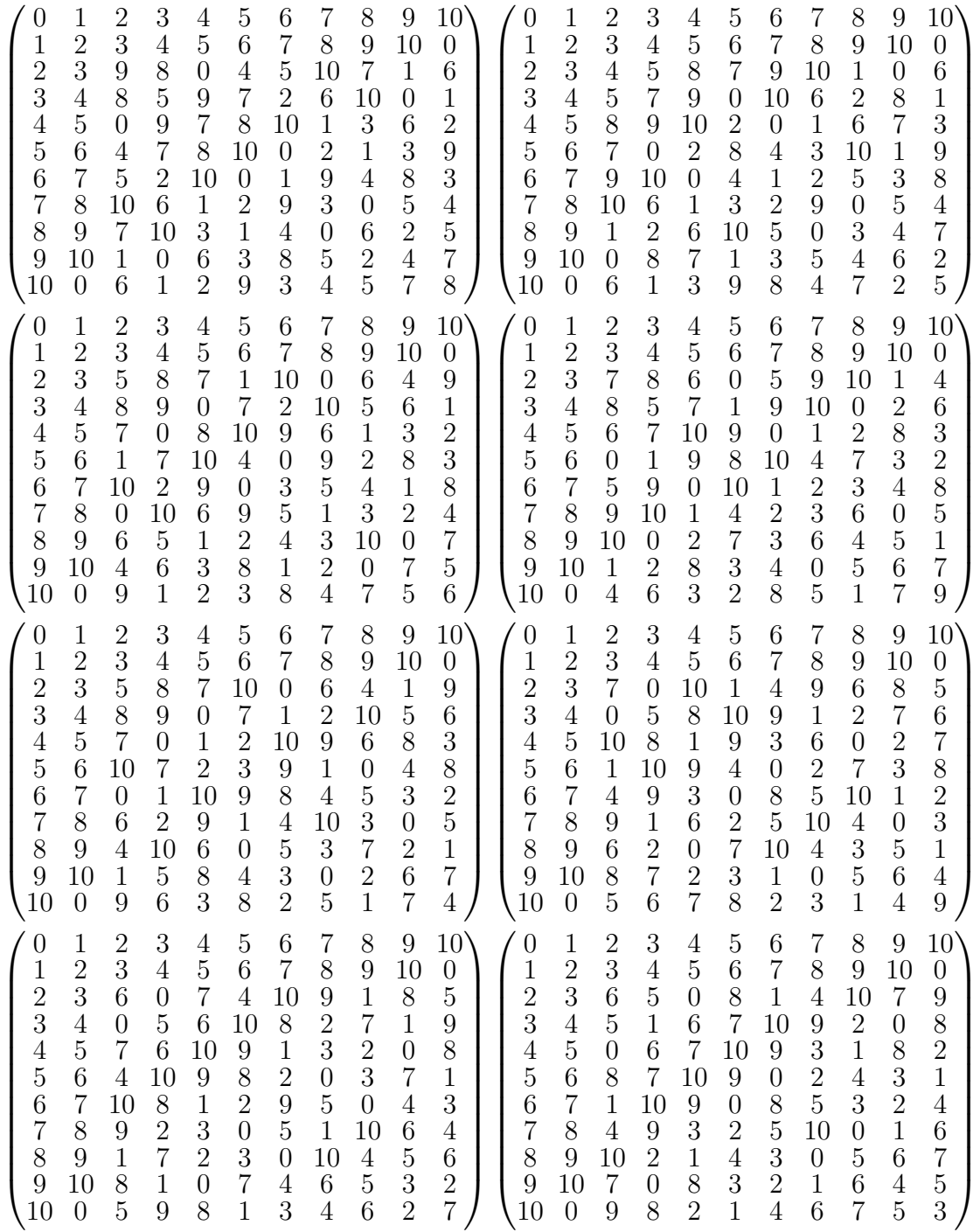

		\[
		\begin{aligned}
			\allowdisplaybreaks
			&\begin{pmatrix}
				0&1&2&3&4&5&6&7&8&9&10 \\
				1&2&3&4&5&6&7&8&9&10&0 \\
				2&3&9&8&0&4&5&10&7&1&6 \\
				3&4&8&5&9&7&2&6&10&0&1 \\
				4&5&0&9&7&8&10&1&3&6&2 \\
				5&6&4&7&8&10&0&2&1&3&9 \\
				6&7&5&2&10&0&1&9&4&8&3 \\
				7&8&10&6&1&2&9&3&0&5&4 \\
				8&9&7&10&3&1&4&0&6&2&5 \\
				9&10&1&0&6&3&8&5&2&4&7 \\
				10&0&6&1&2&9&3&4&5&7&8 \\
			\end{pmatrix}
			\begin{pmatrix}
				0&1&2&3&4&5&6&7&8&9&10\\
				1&2&3&4&5&6&7&8&9&10&0\\
				2&3&4&5&8&7&9&10&1&0&6\\
				3&4&5&7&9&0&10&6&2&8&1\\
				4&5&8&9&10&2&0&1&6&7&3\\
				5&6&7&0&2&8&4&3&10&1&9\\
				6&7&9&10&0&4&1&2&5&3&8\\
				7&8&10&6&1&3&2&9&0&5&4\\
				8&9&1&2&6&10&5&0&3&4&7\\
				9&10&0&8&7&1&3&5&4&6&2\\
				10&0&6&1&3&9&8&4&7&2&5
			\end{pmatrix} \\
			&\begin{pmatrix}
				0&1&2&3&4&5&6&7&8&9&10\\
				1&2&3&4&5&6&7&8&9&10&0\\
				2&3&5&8&7&1&10&0&6&4&9\\
				3&4&8&9&0&7&2&10&5&6&1\\
				4&5&7&0&8&10&9&6&1&3&2\\
				5&6&1&7&10&4&0&9&2&8&3\\
				6&7&10&2&9&0&3&5&4&1&8\\
				7&8&0&10&6&9&5&1&3&2&4\\
				8&9&6&5&1&2&4&3&10&0&7\\
				9&10&4&6&3&8&1&2&0&7&5\\
				10&0&9&1&2&3&8&4&7&5&6
			\end{pmatrix}
			\begin{pmatrix}
				0&1&2&3&4&5&6&7&8&9&10\\
				1&2&3&4&5&6&7&8&9&10&0\\
				2&3&7&8&6&0&5&9&10&1&4\\
				3&4&8&5&7&1&9&10&0&2&6\\
				4&5&6&7&10&9&0&1&2&8&3\\
				5&6&0&1&9&8&10&4&7&3&2\\
				6&7&5&9&0&10&1&2&3&4&8\\
				7&8&9&10&1&4&2&3&6&0&5\\
				8&9&10&0&2&7&3&6&4&5&1\\
				9&10&1&2&8&3&4&0&5&6&7\\
				10&0&4&6&3&2&8&5&1&7&9
			\end{pmatrix} \\
			&\begin{pmatrix}
				0&1&2&3&4&5&6&7&8&9&10\\
				1&2&3&4&5&6&7&8&9&10&0\\
				2&3&5&8&7&10&0&6&4&1&9\\
				3&4&8&9&0&7&1&2&10&5&6\\
				4&5&7&0&1&2&10&9&6&8&3\\
				5&6&10&7&2&3&9&1&0&4&8\\
				6&7&0&1&10&9&8&4&5&3&2\\
				7&8&6&2&9&1&4&10&3&0&5\\
				8&9&4&10&6&0&5&3&7&2&1\\
				9&10&1&5&8&4&3&0&2&6&7\\
				10&0&9&6&3&8&2&5&1&7&4
			\end{pmatrix}
			\begin{pmatrix}
				0&1&2&3&4&5&6&7&8&9&10\\
				1&2&3&4&5&6&7&8&9&10&0\\
				2&3&7&0&10&1&4&9&6&8&5\\
				3&4&0&5&8&10&9&1&2&7&6\\
				4&5&10&8&1&9&3&6&0&2&7\\
				5&6&1&10&9&4&0&2&7&3&8\\
				6&7&4&9&3&0&8&5&10&1&2\\
				7&8&9&1&6&2&5&10&4&0&3\\
				8&9&6&2&0&7&10&4&3&5&1\\
				9&10&8&7&2&3&1&0&5&6&4\\
				10&0&5&6&7&8&2&3&1&4&9
			\end{pmatrix} \\
			&\begin{pmatrix}
				0&1&2&3&4&5&6&7&8&9&10\\
				1&2&3&4&5&6&7&8&9&10&0\\
				2&3&6&0&7&4&10&9&1&8&5\\
				3&4&0&5&6&10&8&2&7&1&9\\
				4&5&7&6&10&9&1&3&2&0&8\\
				5&6&4&10&9&8&2&0&3&7&1\\
				6&7&10&8&1&2&9&5&0&4&3\\
				7&8&9&2&3&0&5&1&10&6&4\\
				8&9&1&7&2&3&0&10&4&5&6\\
				9&10&8&1&0&7&4&6&5&3&2\\
				10&0&5&9&8&1&3&4&6&2&7
			\end{pmatrix}
			\begin{pmatrix}
				0&1&2&3&4&5&6&7&8&9&10\\
				1&2&3&4&5&6&7&8&9&10&0\\
				2&3&6&5&0&8&1&4&10&7&9\\
				3&4&5&1&6&7&10&9&2&0&8\\
				4&5&0&6&7&10&9&3&1&8&2\\
				5&6&8&7&10&9&0&2&4&3&1\\
				6&7&1&10&9&0&8&5&3&2&4\\
				7&8&4&9&3&2&5&10&0&1&6\\
				8&9&10&2&1&4&3&0&5&6&7\\
				9&10&7&0&8&3&2&1&6&4&5\\
				10&0&9&8&2&1&4&6&7&5&3
			\end{pmatrix} 
		\end{aligned}
		\]
		\caption{\label{f:LS45}Eight symmetric row-Hamiltonian Latin squares}
	\end{figure}

	Wanless~\cite{rowhamcount} observed that $11$ is the smallest order for
	which a Latin square with $\nu=4$ exists. \tref{t:LSenum} tells us
	that there are 12 species of Latin squares of order $11$ with $\nu=4$,
	which we now catalogue.  For $m\ge1$ and $b\ge0$ define
	$\Z_{m,b}=\Z_m\cup\{\infty_1,\infty_2,\dots,\infty_b\}$ and
	\[
	z^+ = \begin{cases}
		z+1 & \text{if } z \in \Z_m, \\
		z &  \text{otherwise}.
	\end{cases}
	\]
	A \emph{bordered diagonally cyclic Latin square} (BDCLS) of order
	$m+b$ is a Latin square $L$ of order $m+b$ which satisfies the rule
	that if $(i,j,k)$ is an entry of $L$ then so is $(i^+,j^+,k^+)$. Here
	we are using $\Z_{m,b}$ as the set of row indices, column indices and
	symbols.  If $b=0$ then $L$ is a \emph{diagonally cyclic Latin square} (DCLS).
	For $b\in\{0,1\}$, a BDCLS is uniquely determined by
	its first row \cite{DCLS}. There are four species with $\nu=4$ that
	contain a BDCLS of order $11$. The first row of a BDCLS representative
	for each such species is given below.
	\begin{align}
		&(0, 10, 4, 8, 7, 6, 1, 3, 5, 2, 9), \label{e:fr1}\\
		&(0, 2, 6, 8, 7, \infty_1, 3, 5, 4, 1, 9), \label{e:fr2}\\
		&(0, 3, 7, 9, 8, \infty_1, 4, 6, 5, 2, 1), \text{ and } \label{e:fr3}\\
		&(\infty_1, 1, 9, 7, 5, 3, 8, 6, 4, 2, 0). \label{e:fr4}
	\end{align}
	The DCLS whose first row is \eref{e:fr1} comes from the only
	known infinite family of Latin squares with $\nu=4$
	constructed in~\cite{nu4}. The BDCLS in \eref{e:fr1},
	\eref{e:fr2} and \eref{e:fr3} are each symmetric so, by
	\lref{l:symmnu}, each species gives rise to a single
	isomorphism class of perfect $1$-factorisations.  In contrast,
	the species represented by \eref{e:fr4} gives rise to two
	isomorphism classes of perfect $1$-factorisations.  The
	remaining eight species all contain symmetric Latin
	squares. Figure~\ref{f:LS45} provides a symmetric
	representative of each of these species. As a consequence of
	their symmetry and \lref{l:symmnu}, they also each give rise
	to a single isomorphism class of perfect $1$-factorisations.
	
	The seven species containing atomic Latin squares of order 11
	were catalogued in~\cite{atom11}. From that study it can be
	inferred that they give rise to $12$ isomorphism classes of
	perfect $1$-factorisations. Of course, any species with
	$\nu=2$ can give rise to only a single isomorphism class of
	perfect $1$-factorisations. This accounts for the
	$687\,121=687\,096+13+12$ perfect $1$-factorisations of
	$K_{11,11}$ up to isomorphism. One representative from each
	species containing row-Hamiltonian Latin squares of order 11
	can be found at \cite{wwww}.
	
	Up to paratopism, there are nine row-Hamiltonian Latin squares
	of order $11$ that have trivial autotopism group but
	non-trivial autoparatopism group, and which give rise to
	perfect $1$-factorisations with trivial automorphism
	group. They are the eight squares given in \fref{f:LS45} and a
	symmetric atomic Latin square in the class
	$\mathcal{A}_{11}^5$ from \cite{atom11}. There are two
	isomorphism classes of perfect $1$-factorisations which arise
	from $\mathcal{A}_{11}^5$. One of these has trivial
	automorphism group and the other has automorphism group of
	cardinality $2$.

	We have already given details for the Latin squares reported in
	\tbref{tb:autpargrps} with $\nu>2$. The most symmetric species with
	$\nu=2$ is represented by the DCLS with first
	row
	$ (0, 2, 8, 5, 7, 1, 10, 4, 6, 3, 9)$.
	It has an autotopism that applies the permutation
	$(0,10)(1,9)(2,8)(3,7)(4,6)$ to the rows, columns and symbols.
	Together with the diagonally cyclic symmetry,
	this generates an autotopism group of order $22$.
	The next most symmetric Latin square from \tbref{tb:autpargrps} with $\nu=2$
	is
	\[
	\begin{pmatrix}
		5&7&0&4&9&6&10&8&2&1&3\\
		10&5&8&1&0&7&4&6&9&3&2\\
		1&6&5&9&2&8&3&0&7&10&4\\
		3&2&7&5&10&9&0&4&1&8&6\\
		6&4&3&8&5&10&7&1&0&2&9\\
		4&0&1&2&3&5&8&9&10&6&7\\
		2&9&6&10&4&1&5&3&8&7&0\\
		0&3&10&7&6&2&1&5&4&9&8\\
		7&1&4&6&8&3&9&2&5&0&10\\
		9&8&2&0&7&4&6&10&3&5&1\\
		8&10&9&3&1&0&2&7&6&4&5\\
	\end{pmatrix}.
	\]
	Its autotopism group is isomorphic to the dihedral group of order $10$.

	\section{Invariants}\label{s:invar}
	
	Let $\mathcal{B}$ be the set of row-Hamiltonian Latin squares
	of order $11$, and let $\mathcal{R}(\mathcal{B})$ be the set
	of species representatives of $\mathcal{B}$ that we
	generated. Similarly, let $\mathcal{D}$ be the set of perfect
	$1$-factorisations of $K_{11,11}$, and let
	$\mathcal{R}(\mathcal{D})$ be the set of isomorphism class
	representatives of $\mathcal{D}$ that we generated.  In this
	section we discuss some old and new invariants, and examine
	how useful they are for distinguishing elements of
	$\mathcal{B}$ and elements of $\mathcal{D}$.  A \emph{complete
		species invariant on $\mathcal{B}$} is a function
	$\mathcal{I}$ on $\mathcal{B}$
	such that $\mathcal{I}(L_1)=\mathcal{I}(L_2)$ if
	and only if Latin squares $L_1$ and $L_2$ are paratopic. A
	\emph{complete isomorphism class invariant} for
	$1$-factorisations can be defined similarly.
	
	Let $L$ be a Latin square of order $n$.  A \emph{transversal}
	of $L$ is a selection of $n$ of its entries such that no two
	entries share a row, column or symbol. Let $N(L)$ denote the
	number of transversals of $L$. It is immediate that $N$ is a
	species invariant.
	
	Let $L$ be a Latin square with symbol set $S$ of cardinality
	$n$. Define $G = G(L)$ to be a digraph with vertex set $S^3$
	such that each vertex has a unique outgoing arc. The arc from
	$(a, b, c)$ goes to the triple $(x,y,z)$ where $(a,b,z)$,
	$(a,y,c)$ and $(x,b,c)$ are entries of $L$. The graph $G$ is
	called the \emph{train} of $L$, and the isomorphism class of
	$G$ is a species invariant~\cite{cycatom}. Thus, the indegree
	sequence of $G$ (a sorted list of the indegrees of the
	vertices) is also a species invariant. Denote this indegree
	sequence by $I(L)$.
	
	Recall that a row cycle of a Latin square is a $2 \times k$
	subrectangle that contains no proper subrectangles. We can
	analogously define \emph{column cycles} and \emph{symbol
		cycles}, and taking conjugates interchanges these
	objects. For a Latin square $L$ let $C(L)$ be a sorted list of
	the lengths of its row, column and symbol cycles. Then $C$ is
	a species invariant. Also define $S(L)$ to be a multiset
	consisting of three sorted lists, one giving the lengths of
	its row cycles, one giving the lengths of its column cycles
	and one giving the lengths of its symbol cycles. Then $S$ is
	also a species invariant.
	
	We determined how well the above invariants distinguish
	squares in $\mathcal{B}$ and obtained the following
	results. When applied to every square in
	$\mathcal{R}(\mathcal{B})$:
	\begin{itemize}
		\item $N$ took $630$ values,
		\item $I$ took $283\,518$ values,
		\item $C$ took $151\,412$ values,
		\item $S$ took $675\,110$ values,
		\item $(I, C)$ took $687\,069$ values,
		\item $(N, I, C)$ took $687\,115$ values, thus is a complete invariant on $\mathcal{B}$,
		\item $(I, S)$ took $687\,115$ values, thus is a complete invariant on $\mathcal{B}$.
	\end{itemize} 
	
	Let $\F$ and $\E$ be non-isomorphic perfect $1$-factorisations of $K_{11,11}$ such that $\L(\F)$ is paratopic to $\L(\E)$. Since each of $N$, $I$, $C$ and $S$ are well known species invariants, they cannot possibly distinguish between $\F$ and $\E$. So we now define a new invariant, which is useful for distinguishing such perfect $1$-factorisations.
	
	Let $\F = \{f_1, f_2, \ldots, f_n\}$ be a perfect $1$-factorisation of $K_{n, n}$. Let $\{i, j, k\} \subseteq \{1, 2, \ldots, n\}$ with $i < j$ and $k \notin \{i, j\}$. Let $\F_{i, j}$ denote the subgraph of $K_{n, n}$ with edge set $f_i \cup f_j$. Since $\F$ is perfect, $\F_{i, j}$ forms a Hamiltonian cycle in $K_{n, n}$. For each edge $e \in f_k$, define $p_{i, j, k, e}$ to be the distance between the endpoints of $e$ in $\F_{i, j}$. Define
	\[
	P(\F) = \sum_{i < j} \sum_{k \neq i, j} \prod_{e \in f_k} p_{i, j, k, e}.
	\]
	Then $P$ is invariant on isomorphism classes of perfect $1$-factorisations of $K_{n, n}$.
	
	When applied to every element of $\mathcal{R}(\mathcal{D})$,
	$P$ took $687\,115$ values. The six pairs of elements in
	$\mathcal{R}(\mathcal{D})$ on which $P$ coincided can
	be found at \cite{wwww}. For any invariant $\mathcal{I} \in
	\{N,C,I,S\}$, the pair $(P, \mathcal{I})$ took $687\,121$
	values, thus formed a complete invariant on $\mathcal{D}$.

	\section{The classical families of perfect $1$-factorisations}\label{s:GA}

	The main purpose of this section is to revisit the classical
	families of perfect $1$-factorisations to tie up a loose end
	in the literature, which we do in \tref{t:loosend} below.  As
	mentioned in~\sref{s:intro}, given a perfect $1$-factorisation
	of the complete graph $K_{n+1}$, there is a known method for
	constructing perfect $1$-factorisations of $K_{n,n}$. We refer to this as
	Kotzig's construction. Let $n$ be an odd integer, let $V$ be
	the vertex set of $K_{n+1}$, and suppose that $\F$ is a
	perfect $1$-factorisation of $K_{n+1}$. For distinct $x$ and
	$y$ in $V$ let $h_{x, y}$ denote the unique $1$-factor in $\F$
	containing the edge $xy$. Fix a vertex $v \in V$ called the
	\emph{root} vertex. We associate to the pair $(\F, v)$ a Latin
	square of order $n$, denoted by $\L(\F, v)$, whose row index
	set, column index set and symbol set is $V \setminus \{v\}$,
	and is defined by
	\[
	\L(\F, v)_{i, j} = \begin{cases}
		i & \text{if } j = i,\\
		k & \text{if } j \neq i, \text{ where } k \in V \setminus \{v\} \text{ is such that }kv\in h_{i,j}.
	\end{cases}
	\]
	Then $\L(\F, v)$ is a symmetric Latin square whose
	$(321)$-conjugate is row-Hamiltonian and hence encodes a
	perfect $1$-factorisation of $K_{n,n}$. \lref{l:symmnu} implies
	that $\nu(\L(\F, v)) \in \{2,6\}$. Furthermore, if
	$\{u,v\}\subseteq V$ then $\L(\F,v)$ is paratopic to
	$\L(\F,u)$ if and only if there is an automorphism of $\F$
	that maps $v$ to $u$. See~\cite{symm1facs} for more details.
	
	We now discuss the known infinite families of row-Hamiltonian Latin
	squares that come from the construction given above. For each prime $p
	\geq 11$ there are two known non-isomorphic perfect $1$-factorisations
	of $K_{p+1}$ which come from infinite families. One is due to Kotzig
	and is commonly denoted by $GK_{p+1}$. The other is due to Bryant,
	Maenhaut, and Wanless~\cite{newatom}, which we will denote by
	$GB_{p+1}$. There are two species that contain Latin squares
	$\L(GK_{p+1}, v)$ for some root vertex $v$, and there are three other
	species that contain Latin squares $\L(GB_{p+1}, v)$ for some root
	vertex $v$. If $2$ is primitive modulo $p$ then all five of these
	species have $\nu=6$. If $2$ is not primitive modulo $p$ then two of
	these species have $\nu=6$ and the remaining three have $\nu=2$,
	see~\cite{newatom, M3atom, rowhamcount}. There is a well known perfect
	$1$-factorisation of $K_{2p}$ for every odd prime $p$, commonly
	denoted by $GA_{2p}$. Kotzig~\cite{Kotzconj} stated that $GA_{2p}$ is
	perfect for every odd prime $p$, and a proof was provided by
	Anderson~\cite{GAperfect}. For each odd prime $p$, every Latin square
	of the form $\L(GA_{2p}, v)$ lies in the same species. Our goal for
	this section is to show that this species has $\nu=2$ unless $p=3$, in
	which case it has $\nu=6$.
	
	There are some infinite families of row-Hamiltonian Latin
	squares that do not come from perfect $1$-factorisations of
	complete graphs.  For each prime $p \geq 11$, Bryant, Maenhaut
	and Wanless \cite{p1fbip} constructed $(p-1)/2$ species of
	order $p^2$ with $\nu=2$. Allsop and Wanless \cite{nu4}
	constructed, for each prime $p\not\in\{3,19\}$ with $p \equiv
	1 \bmod 8$ or $p \equiv 3 \bmod 8$, a Latin square of order
	$p$ with $\nu=4$. There are also some sporadic examples of
	row-Hamiltonian Latin squares~\cite{P1F16,cycatom}.

	We now return to the family $GA_{2p}$. Let $p$ be an odd prime
	and let the vertex set of $K_{2p}$ be
	$\Z_{p}\times\{1,2\}$. For $i\in\Z_p$ define
	\[
	f_i = \{(i+j,1)(i-j,1),(i+j,2)(i-j,2):j\in\{1,2,\dots,(p-1)/2\}\}\cup\{(i,1)(i,2)\}.
	\]
	For $i \in \Z_p \setminus \{0\}$ define
	\[
	g_i = \{(j,1)(i+j,2):j\in\Z_p\}.
	\]
	Then
	\[
	GA_{2p} = \{f_i : i \in \Z_p\} \cup \{g_i : i \in \Z_p \setminus \{0\}\}
	\]
	is a perfect $1$-factorisation of $K_{2p}$. 
	
	Anderson~\cite{GAaut} showed that the automorphism group of $GA_{2p}$
	acts transitively on the vertices of $K_{2p}$. Since we are only
	interested in the species of Latin square obtained from $GA_{2p}$, we
	may decide to work with the root vertex $v=(-1,2)$. Define
	$\LL_p=\L(GA_{2p},v)$. We can give a more explicit
	definition of $\LL_p$.
	
	\begin{lem}\label{l:Lp}
		The square $\LL=\LL_p$ is defined by
		\[
		\LL_{(x, z), (y, w)}=\begin{cases}
			(x, z) & \text{if } (x, z)=(y, w),\\
			(x+y+1,2) & \text{if } z=w \text{ and } x+y+2 \neq 0,\\
			(-1,1) & \text{if } z=w \text{ and } x+y+2=0,\\
			(2x+1,2) & \text{if } z \neq w, \text{ and } x=y,\\
			(x-y-1,1) & \text{if } z=1, w=2, \text{ and } x \neq y,\\
			(y-x-1,1) & \text{if } z=2, w=1 \text{ and } x \neq y.
		\end{cases}
		\]
	\end{lem}
	
	\begin{proof}
		Let $((x, z), (y, w)) \in (\Z_p \times \{1,2\})^2$
		with $(x, z) \neq (y, w)$. First suppose that
		$z=1=w$. Let $i=2^{-1}(x+y)\in\Z_p$ and note that
		$(x,z)(y,w)\in f_i$.  If $x+y+2=0$ then $i=-1$ and so
		$(-1,1)(-1,2)\in f_i$. Hence
		$\LL_{(x,z),(y,w)}=(-1,1)$. Now suppose that $x+y+2\neq 0$.
		Let $j = i+1 \in \Z_p$ so that $i-j = -1$. Then
		$i+j=2i+1=x+y+1$ and thus $(x+y+1,2)(-1,2)\in f_i$.
		Hence $\LL_{(x,z),(y,w)}=(x+y+1,2)$. Similar
		arguments can be used to prove that the claimed value
		of $\LL_{(x,z),(y,w)}$ is correct when $z=2=w$.
		
		Now assume that $z=1$ and $w=2$. We must distinguish
		two cases depending on whether or not $x=y$. First
		suppose that $x \neq y$. Let $i = y-x$ and note that
		$(x,z)(y,w)\in g_i$.  Setting $i+j=-1$ we obtain
		$j=-i-1=x-y-1$. So $(x-y-1,1)(-1,2)\in g_i$ and thus
		$\LL_{(x,z),(y,w)}=(x-y-1,1)$. Now suppose that
		$y=x$. Then $(x,z)(y,w)\in f_x$. Setting $x-j=-1$
		yields $j=x+1$. Thus $(2x+1,2)(-1,2)\in f_x$ and so
		$\LL_{(x,z),(y,w)}=(2x+1,2)$. Similar arguments can be
		used to prove that the claimed value of
		$\LL_{(x,z),(y,w)}$ is correct when $z=2$ and $w=1$.
	\end{proof}
	
	We are now ready to determine $\nu(\LL_p)$ for each odd prime $p$.
	
	\begin{thm}\label{t:loosend}
		If $p=3$ then $\LL_p$ is atomic and otherwise $\nu(\LL_p)=2$. 
	\end{thm}
	
	\begin{proof}
		By \tbref{tb:smallnLSrowhamatomic}, we know that any Latin square of order $5$ that has $\nu>0$ is atomic, so the theorem is true when $p=3$. Now assume that $p \geq 5$. Using \lref{l:Lp} it is easy to verify that the following ten triples are entries of $\LL_p$:
		\[
		\begin{array}{llll}			
			((0,1),(0,1),(0,1)),
			&&&((0,2),(0,1),(1,2)),\\
			((0,1),(0,2),(1,2)),
			&&&((0,2),(0,2),(0,2)),\\
			((0,1),(-1,1),(0,2)),
			&&&((0,2),(-1,1),(-2,1)),\\
			((0,1),(1,2),(-2,1)),
			&&&((0,2),(1,2),(2,2)),\\
			((0,1),(1,1),(2,2)),
			&&&((0,2),(1,1),(0,1)).			
		\end{array}
		\]
		These entries form a row cycle of length $5$ in
		$\LL_p$ and so $\LL_p$ is not atomic. Since
		$\nu(\LL_p) \in \{2, 6\}$ by \lref{l:symmnu}, this
		proves the lemma.
	\end{proof}
	
	We end by discussing the result of applying Kotzig's
	construction to the perfect $1$-factorisations of complete
	graphs of small order, where we have exhaustive
	catalogues. Suppose that $L=\L(\F,v)$ where $n\le15$ and $\F$
	is any of the perfect $1$-factorisations of $K_{n+1}$, with
	$v$ being one of the vertices of $K_{n+1}$.  By
	\lref{l:symmnu}, we know that $\nu(L)=2$ unless $L$ is atomic.
	We can infer from \cite{rowhamcount} and \cite{P1F16}
	respectively that $\nu(L)=2$ if $n\in\{9,15\}$ since there are
	no symmetric atomic Latin squares of these orders. If $n=13$ there are five species of atomic Latin squares that can be derived from perfect
	$1$-factorisations of $K_{14}$, and these are described
	in \cite{newatom}. The situation for $n=11$ is covered in
	detail in \cite{atom11}; there are six species of atomic Latin squares of order $11$ that come from
	perfect $1$-factorisations of $K_{12}$.  For odd $n\le7$ it is
	known that there is a unique species of atomic Latin squares,
	which can be obtained from $GK_{n+1}$, the unique perfect
	$1$-factorisation of $K_{n+1}$ (although, for $n=7$ it is
	important to choose the root vertex $v$ to be the unique
	vertex that is fixed by all automorphisms of $GK_8$).
	
	\subsection*{Acknowledgements}
	The authors are grateful for access to the MonARCH HPC
	Cluster, where they did their computations. This research was
	supported by Australian Research Council grant DP250101611.


\begin{thebibliography}{20}
		
		\bibitem{Ninf}
		J. Allsop and I. M. Wanless, ``Latin squares without proper subsquares", \textit{arXiv:2310.01923} (2023).
		
		\bibitem{nu4}
		J. Allsop and I. M. Wanless, ``Row-Hamiltonian Latin squares and Falconer varieties", \textit{Proc. Lond. Math. Soc. (3)} 128.1 (2024), Paper No. e12575.
		
		\bibitem{GAperfect}
		B. A. Anderson, ``Finite topologies and Hamiltonian paths", \textit{J. Combin. Theory Ser. B} 14 (1973) 87--93.
		
		\bibitem{GAaut}
		B. A. Anderson, ``Symmetry groups of some perfect 1-factorizations of complete graphs", \textit{Discrete Math.} 18.3 (1977), 227--234.
		
		\bibitem{newatom}
		D. Bryant, B. Maenhaut, and I. M. Wanless, ``New families of atomic Latin squares and perfect 1-factorisations", \textit{J. Combin. Theory Ser. A} 113.4 (2006), 608--624.
		
		\bibitem{p1fbip}
		D. Bryant, B. M. Maenhaut, and I. M. Wanless, ``A family of perfect factorisations of complete bipartite graphs", \textit{J. Combin. Theory Ser. A} 98.2 (2002), 328--342.
		
		\bibitem{P1F14}
		J. H. Dinitz and D. K. Garnick, ``There are 23 nonisomorphic perfect one-factorizations of $K_{14}$", \textit{J. Combin. Des.} 4.1 (1996), 1--4.
		
		\bibitem{P1F10}
		E. N. Gelling and R. E. Odeh, ``On 1-factorizations of the complete graph and the relationship to round robin schedules", \textit{Congr. Numer.} 9 (1974) 213--221.
		
		\bibitem{P1F16}
		M. J. Gill and I. M. Wanless, ``Perfect $1$-factorisations of $K_{16}$", \textit{Bull. Aust. Math. Soc.} 101.2 (2020), 177-185.
		
		\bibitem{Kotzconj}
		A. Kotzig, ``Hamilton graphs and Hamilton circuits", \textit{Theory of Graphs and its Applications (Proc. Sympos. Smolenice, 1963)}, Publ. House Czechoslovak Acad. Sci., Prague, 1964, 63--82.
		
		\bibitem{atom11}
		B. M. Maenhaut and I. M. Wanless, ``Atomic Latin squares of order eleven", \textit{J. Combin. Des.} 12.1 (2004) 12--34.
		
		\bibitem{nauty}
		B. D. McKay and A. Piperno, ``Practical graph isomorphism, {II}", \textit{J. Symbolic Comput.} 60 (2014), 94--112.
		
		\bibitem{P1F162}
		M. Meszka, ``There are 3155 nonisomorphic perfect one-factorizations of $K_{16}$", \textit{J. Combin. Des.} 28.1 (2020), 85--94.
		
		\bibitem{M3atom}
		P. J. Owens and D. A. Preece, ``Some new non-cyclic Latin squares that have cyclic and Youden properties", \textit{Ars Combin.} 44 (1996), 137--148.
		
		\bibitem{P1F12}
		A. P. Petrenyuk and A. Ya. Petrenyuk, ``Intersection of perfect 1-factorizations of complete graphs", \textit{Cybernetics} 16.1 (1980), 6--9.
		
		\bibitem{cycatom}
		I. M. Wanless, ``Atomic Latin squares based on cyclotomic orthomorphisms", \textit{Electron. J. Combin.} 12 (2005), Paper No. 22.
		
		\bibitem{DCLS}
		I. M. Wanless, ``Diagonally cyclic Latin squares", \textit{European J. Combin.} 25.3 (2004), 393--413.
		
		\bibitem{rowhamcount}
		I. M. Wanless, ``Perfect factorisations of bipartite graphs and Latin squares without proper subrectangles", \textit{Electron. J. Combin.} 6 (1999), Paper No. 9.
		
		\bibitem{wwww}
		I. M. Wanless, \textit{User homepage}, \url{https://users.monash.edu.au/~ iwanless/data/}.
		
		\bibitem{symm1facs}
		I. M. Wanless and E. C. Ihrig, ``Symmetries that Latin squares inherit from 1-factorizations", \textit{J. Combin. Des.} 13.3 (2005), 157--172.
		
	\end{thebibliography}
\end{document}